\documentclass{article}
\usepackage{amsmath}
\usepackage{amssymb}
\usepackage{indentfirst}
\usepackage{enumerate}

\newtheorem{lemma}{Lemma}[section]

\newenvironment{proof.}[1][Proof]{\begin{trivlist}
\item[\hskip \labelsep {\bfseries #1}]}{\end{trivlist}}
\begin{document}

\title{The Crucial Constants in the Exponential-type Error Estimates for Multiquadric Interpolations}         
\author{LIN-TIAN LUH}        

\maketitle     

{\bf Abstract}
In this article we explore the exponential-type error bound for
multiquadric and inverse multiquadric interpolations, which was put forward by
Madych and Nelson in 1992. It is of the form $\left\vert
f(x)-s(x)\right\vert \leq \lambda ^{\frac{1}{d}}\left\Vert f\right\Vert _{h},$%
 \ where $0<\lambda <1$ is a constant, $d$ is the fill
distance which roughly speaking measures the spacing of the data
points, $s(x)$ is the interpolating function of $f(x),$ and $h$
denotes the multiquadric or inverse multiquadric. The error bound
converges very fast as $d\rightarrow 0. $ The constant $\lambda $
is very sensitive. A slight change of it will result in a huge
change of the error bound. Unfortunately, $\lambda $ can not be
calculated, or even approximated. In Wendland's book $\left[ 8\right]$ there is also an exponential-type error bound for multiquadric interpolation. However, the crucial constant $\lambda$ there remains unknown. This is a famous question in the
theory of radial basis functions. The purpose of this paper is to
answer the question. Incidentally, we find that the constant $\lambda$ greatly depends on the shape parameter c contained in the multiquadrics. It ushers in a useful criterion for the optimal choice of the shape parameter.  \\
\\
\maketitle {\bf Key words:} radial basis
function, conditionally positive definite function, interpolation,
multiquadric, inverse multiquadric\\
\\
\maketitle {\bf AMS subject classification}: 41A05, 41A15, 41A25,
41A30, 41A63

\begin{center}
\section{Introduction}
\end{center}

Let $h$ be a continuous function on $\mathbb{R}^n$ which is
conditionally
positive definite of order $m$. Given data $\left( x_{j},y_{j}\right) ,$ $%
j=1,...,~N,$ where $X=\left\{ x_{1},...,x_{N}~\right\} $ is a
subset of points in ${\mathbb{R}}^n$ and the $y_{j}$
are real or complex
numbers, the so-called $h$ spline interpolant of these data is the function $%
s$ defined by
\begin{align}
s(x)=p(x)+\sum\limits_{j=1}^{N}c_{j}h(x-x_{j}), \tag{1.1}
\end{align}
where $p(x)$ is a polynomial in $\mathcal{P}_{m-1}$ and the $c_{j}$ are chosen so that
\begin{align}
\sum_{j=1}^{N}c_{j}q(x_{j})=0 \tag{1.2}
\end{align}
for all polynomials $q$ in $\mathcal{P}_{m-1}$ and
\begin{align}
p(x_{i})+ \sum_{j=1}^{N}c_{j}h(x_{i}-x_{j})=y_{i},\quad i=1,...,~N.\tag{1.3}
\end{align}
Here $\mathcal{P}_{m-1}$ denotes the class of those polynomials of \
 ${\mathbb{R}}^n$ of degree $\leq m-1$.

  As pointed out in [4], the linear system of equations (1.2)and (1.3) has a unique solution when $X$ is a determining set for $\mathcal{P}_{m-1}$ and $h$ is strictly conditionally positive definite. For
more details, please see [4]. Therefore, the interpolant $s(x)$ here is well defined because we only deal with multiquadrics and inverse multiquadrics which are conditionally positive definite of order $m\geq 0$.

 A set $X$ is said to be a determining set for $\mathcal{P}_{m-1}$ if $p$ is in $\mathcal{P}_{m-1}$ and $p$ vanishes on $X$
imply that $p$ is identically zero. This kind of sets is also called polynomial-nondegenerate by some authors.

 In this paper $h$ is defined by the formula
\begin{align}
h(x):=\Gamma (-\frac{\beta }{2})(c^{2}+\left\vert x\right\vert^{2})^{\frac{\beta }{2}},\quad \beta \in R\setminus 2N_{\geq 0},\quad c>
0,\tag{1.4}
\end{align}
where $\left\vert x\right\vert $ is the Euclidean norm of $x,$ $\Gamma $ is
\ the classical gamma function and $\beta ,$ $c$ are constants. The function
$h$ is called multiquadric or inverse multiquadric, respectively, depending
on $\beta >0,~$or $\beta <0$.

 In $\left[ 5\right] $ Madych and Nelson obtain bounds on the pointwise
difference between a function $f$ and the $h$ spline which agrees with $f$
on a subset $X$ of $R^{n}$. These estimates involve a parameter that
measures the spacing of the points in $X$ and are $O\left( d^{\mathcal{\ell
}}\right) $ as $d\rightarrow 0$ where $l$ depends on $h$. Later in 
$\left[ 6\right] $ they find that for multiquadrics and inverse multiquadrics, the
estimate can be improved to $O\left( \lambda ^{\frac{1}{d}}\right) $ as $d\rightarrow 0,$ where $\lambda $ is a constant which satisfies $0<\lambda
<1.$ The conditions on $f$ are the same as those in $\left[ 5\right]$.

 There is something we have to clarify. Because this area of material is quite complicated and not easy to understand, we try to maintain the features and notations of [6] in order to avoid unnecessary troubles. There are two advantages for doing so. First, it will be easier to find the differences between the results of [6] and this paper. Second, it will help the reader understand Madych and Nelson's great works via this paper. For the first point, we emphasize on the constant $\rho$ which also appears in [6]. The calculation of $\rho$ is a big problem and has been considered to be a hard question. Madych and Nelson only point out that such a constant exists, but no tool is offered to calculate it. We solve it in this paper. Once this is solved, the calculations of other constants, especially $\lambda$, can be easily achieved by following the route of [6]. Nothing has to be changed; otherwise it will become messy. Consequently the reader shouldn't be surprised if they find too many counterparts in the two papers.\\

\subsection{An Inequality for Multivariate Polynomials}

 Our theory involves an inequality for multivariate polynomials whose proof is very deep and is based on the theory of algebraic geometry. In order to avoid digression, we will not talk too much about its construction. We cite it directly
from $\left[ 6\right] $ and omit its proof.

\begin{lemma}
For $n=1,2,...,~$define $\gamma _{n}$ by the formulas $\gamma
_{1}=2 $ and, if $n>1,$ $\gamma _{n}=2n(1+\gamma _{n-1}).$ Let $Q$
be a cube in ${\mathbb{R}}^n$ that is subdivided into $q^{n}$
identical subcubes. Let $Y$ be a set of $q^{n}$ points obtained by
selecting a point from each of those subcubes. If $q\geq \gamma
_{n}(k+1),$ then for all $p$ in $\mathcal{P}_{k}$
\[
\sup_{x \in Q}\left\vert p(x)\right\vert\leq e^{2n \gamma_{n}(k+1)}\sup_{y \in Y}\left\vert p(y)\right\vert.
\]
\end{lemma}

\subsection{Function Space and Interpolation Setting}

 Our calculation of $\lambda$ involves the basic theory of interpolation introduced by Madych and Nelson in [5]. In order to make this paper more readable, let's review some basic ingredients and notations.

 The space of complex-valued functions on ${\mathbb{R}}^n$
that are compactly supported and infinitely differentiable is
denoted by $\mathcal{D}$. The Fourier transform of a function
$\phi $ in $\mathcal{D}$ is
\[
\hat{\phi}(\xi )=\int e^{-i<x,\xi >}\phi(x)dx.
\]
A continuous function $h$ is conditionally positive definite of order $m$ if
\[
\int h(x)\phi (x)\ast \tilde{\phi}(x)dx\geq 0
\]
holds whenever $\phi =p(D)\psi $ with $\psi $ in $\mathcal{D}$ and $p(D)$ a
linear homogeneous constant coefficient differential operator of order $m$.
Here $\tilde{\phi}(x)=\overline{\phi (-x)}$ and $\ast$ denotes the
convolution product
\[
\phi _{1}\ast \phi _{2}(t)=\int \phi _{1}(x)\phi _{2}(t-x)dx.
\]
As pointed out in $\left[ 5\right]$ , this definition of conditional
positive definiteness is equivalent to that of $\left[ 4\right], $ which is
generally used.

 If $h$ is a continuous and conditionally positive definite
function of order $m$, the Fourier transform of $h$ uniquely determines a
positive Borel measure $\mu$ on $R^{n}$ $\backslash \left\{ 0\right\}$ and constants
$a_{r}$,$\left\vert r\right\vert$ =2$m$ as follows: For all $\psi \in \mathcal{D}$
\begin{align}
\int h(x)\psi (x)dx = \int \left\{ \hat{\psi}(\xi )-\hat{\chi}(\xi )\sum_{\left\vert \gamma \right\vert <2m}D^{\gamma}\hat{\psi}(0)\frac{\xi ^{\gamma}}{\gamma!}\right\} d\mu (\xi )\tag{1.5}
\end{align}
\[
+\sum_{\left\vert \gamma\right\vert \leq 2m}D^{\gamma}\hat{\psi}(0)\frac{a_{\gamma}}{\gamma!},
\]
where for every choice of complex numbers $c_{\alpha },\left\vert \alpha
\right\vert =m$,
\[
\sum_{\left\vert \alpha \right\vert =m}\sum_{\left\vert
\beta \right\vert =m}a_{\alpha +\beta }c_{\alpha }\overline{c_{\beta }}\geq 0.
\]
Here $\chi $ is a function in $\mathcal{D}$ such that $1-\hat{\chi}(\xi )$
has a zero of order $2m+1$ at $\xi =0;$ both of the integral $%
\int_{0<\left\vert \xi \right\vert <1}\left\vert \xi \right\vert ^{2m}d\mu
(\xi ),$ $\int_{\left\vert \xi \right\vert \geq 1}d\mu (\xi )$ are finite.
The choice of $\chi $ affects the value of the coefficients $a_{\gamma }$
for $\left\vert \gamma \right\vert <2m.$

 In this paper, the space of the interpolated functions is denoted by $\mathcal{C}_{h,m}$ which some people call the native space. If
\[
\mathcal{D}_{m}=\left\{ \phi \in \mathcal{D}:\int x^{\alpha }\phi
(x)dx=0\quad for~all~\left\vert \alpha \right\vert <m\right\} ,
\]
then {\bf $\mathcal{C}_{h,m}$} is the class of those continuous functions $f$
which satisfy
\begin{align}
\left\vert \int f(x)\phi (x)dx\right\vert \leq
c(f)\left\{ \int h(x-y)\phi (x)\overline{\phi (y)}dxdy\right\} ^{\frac{1}{2}}\tag{1.6}
\end{align}
for some constant $c(f)$ and all $\phi $ in $\mathcal{D}_{m}.$ If $f\in
\mathcal{C}_{h,m}$, let $\left\Vert f\right\Vert _{h}$ denote the smallest
constant $c(f)$ for which $\left( 1.6\right) $ is true. Recall that $%
\left\Vert f\right\Vert _{h}$ is a semi-norm and $\mathcal{C}_{h,m}$ is a
semi-Hilbert space; in the case $m=0$ it is a norm and a Hilbert space
respectively. The characterizations of the native space can be found in [1],[2],[3],[4],[5] and [8]. 

\begin{center}
\section{Main Results}
\end{center}

We first recall that the function $h$ defined in $\left( 1.4\right) $ is
conditionally positive definite of order $m=0$ if $\beta <0,$ and $%
m=\left\lceil \frac{\beta }{2}\right\rceil $ if $\beta >0.$ This can be
found in $\left[ 7\right] $ and many relevant papers. Then we have the
following lemma.\\

\begin{lemma}
Let $h$ be as in $\left( 1.4\right) $ and $m$ be its order of
conditional positive definiteness. There exists a positive constant $\rho $
such that
\begin{align}
\int_{R^{n}}\left\vert \xi \right\vert ^{k}d\mu \left( \xi
\right) \leq \left( \sqrt{2}\right) ^{n+\beta +1}\cdot \left( \sqrt{\pi }%
\right) ^{n+1}\cdot n\alpha _{n}\cdot c^{\beta -k}\cdot \triangle _{0}\cdot
\rho ^{k}\cdot k!\tag{2.1}
\end{align}
for all integer $k\geq 2m+2$ where $\mu $ was defined in $\left( 1.5\right) ,$
$\alpha _{n}$ denotes the volume of the unit ball in $R^{n},$ c is as in $%
\left( 1.4\right) ,$ and $\triangle _{0}$ is a positive constant.
\end{lemma}
\em \bf Proof.\rm \quad
Let $K_{\nu }$ denote the modified Bessel function of the second kind. Then
\begin{eqnarray}
& & \int_{R^{n}}\left\vert \xi \right\vert ^{k}d\mu \left( \xi \right) \nonumber \\
& =& \int_{R^{n}}\left\vert \xi \right\vert ^{k}\cdot 2\pi
^{\frac{n}{2}}\cdot
\left( \frac{\left\vert \xi \right\vert }{2c}\right) ^{-\frac{n+\beta }{2}%
}\cdot K_{\frac{n+\beta }{2}}\left( c\left\vert \xi \right\vert \right) d\xi \nonumber \\
& = & 2\pi ^{\frac{n}{2}}\left( \frac{1}{2c}\right)
^{-\frac{n+\beta }{2}}\cdot
\int_{R^{n}}\left\vert \xi \right\vert ^{k-\frac{n+\beta }{2}}\cdot K_{\frac{%
n+\beta }{2}}\left( c\left\vert \xi \right\vert \right) d\xi \nonumber \\
& \sim & \frac{\sqrt{\pi }}{\sqrt{2}}\cdot 2\pi ^{\frac{n
}{2}}\cdot \left( \frac{1}{2c}\right) ^{-\frac{n+\beta
}{2}}\int_{R^{n}}\left\vert \xi \right\vert ^{k-\frac{n+\beta
}{2}}\cdot \frac{1}{\sqrt{c\left\vert \xi
\right\vert }\cdot e^{c\left\vert \xi \right\vert }}d\xi \nonumber \ (See\ p.11)\\
& = & \frac{\sqrt{\pi }}{\sqrt{2}}\cdot 2\pi ^{\frac{n}{2}}\cdot \left( \frac{1}{%
2c}\right) ^{-\frac{n+\beta }{2}}\cdot n\cdot \alpha _{n}\int_{_{0}}^{\infty
}r^{k-\frac{n+\beta }{2}}\cdot \frac{r^{n-1}}{\sqrt{c\left\vert r
\right\vert }\cdot e^{c\left\vert r \right\vert }}dr \nonumber \\
& = & \frac{\sqrt{\pi }}{\sqrt{2}}\cdot 2\pi ^{\frac{n}{2}}\cdot
\left( 2c\right) ^{\frac{n+\beta }{2}}\cdot
\frac{1}{\sqrt{c}}\cdot n\cdot \alpha
_{n}\int_{_{0}}^{\infty }\frac{r^{k+\frac{n-\beta -3}{2}}}{e^{cr}}dr \nonumber \\
& = & \frac{\sqrt{\pi }}{\sqrt{2}}\cdot 2\pi ^{\frac{n}{2}}\cdot
\left( 2c\right) ^{\frac{n+\beta }{2}}\cdot
\frac{1}{\sqrt{c}}\cdot n\cdot \alpha
_{n}\cdot \frac{1}{c^{k+\frac{n-\beta -1}{2}}}\int_{_{0}}^{\infty }\frac{r^{k+\frac{n-\beta -3}{2}}}{e^{r}}dr \nonumber \\
& = & 2^{\frac{n+\beta +1}{2}}\cdot \pi ^{\frac{n+1}{2}}\cdot
n\cdot \alpha_{n}\cdot c^{\beta -k}\int_{_{0}}^{\infty
}\frac{r^{k^{^{\prime }}}}{e^{r}}dr \ where \quad k^{^{\prime
}}=k+\frac{n-\beta -3}{2}.\nonumber
\end{eqnarray}
Note that if $\beta <0,$ then $m=0$ and $k\geq 2m+2=2$. This implies $%
k^{^{\prime }}>0.$ If $\beta >0,$ then $m=\left\lceil \frac{\beta }{2}%
\right\rceil $ and $k\geq 2m+2=2\left\lceil \frac{\beta }{2}\right\rceil +2.$
This implies $k^{^{\prime }}>0.$ In any case $k^{^{\prime }}>0.$

 Now we divide the proof into three cases. Let $k^{^{\prime \prime
}}=\left\lceil k^{^{\prime }}\right\rceil $ which is the smallest integer
greater than or equal to $k^{^{\prime }}.$\\
\\
\em \bf Case1. \rm
Assume $k^{^{\prime \prime }}>k.$ Let $k^{^{\prime \prime}}=k+s.$ Then
\[
\int_{0}^{\infty }\frac{r^{k^{^{\prime }}}}{e^{r}}dr\leq \int_{0}^{\infty }%
\frac{r^{k^{^{\prime \prime }}}}{e^{r}}dr=k^{^{\prime \prime
}}!=(k+s)(k+s-1)\cdot \cdot \cdot (k+1)k!
\]
and
\[
\int_{0}^{\infty }\frac{r^{k^{^{\prime }}+1}}{e^{r}}dr\leq \int_{0}^{\infty
}\frac{r^{k^{^{\prime \prime }}+1}}{e^{r}}dr=(k^{^{\prime \prime
}}+1)!=(k+s+1)(k+s)\cdot \cdot \cdot (k+2)(k+1)k!.
\]
Note that
\[
\frac{(k+s+1)(k+s)\cdot \cdot \cdot (k+2)}{(k+s)(k+s-1)\cdot \cdot \cdot
(k+1)}=\frac{k+s+1}{k+1}.
\]
\rm{(i)}\em Assume $\beta <0$. Then $m=0$ and $k\geq 2$. This
gives
\[
\frac{k+s+1}{k+1}\leq \frac{3+s}{3}.
\]
Let\quad $\rho =\frac{3+s}{3}$. Then
\[
\int_{0}^{\infty }\frac{r^{k^{^{\prime
\prime }}+1}}{e^{r}}dr\leq \triangle _{0}\cdot \rho ^{k+1}\cdot (k+1)!.
\]
if $\int_{0}^{\infty }\frac{r^{k^{^{\prime \prime }}}}{e^{r}}dr\leq
\triangle _{0}\cdot \rho ^{k}\cdot k!$. The smallest $k^{^{\prime \prime }}$ is $k_{0}^{^{\prime \prime }}=2+s$. Now,
\begin{eqnarray*}
\quad \int_{0}^{\infty }\frac{r^{k_{0}^{\prime \prime }}}{e^{r}}dr
& = & k_{0}^{^{\prime\prime }}!=(2+s)(2+s-1)\cdot \cdot \cdot (3)\cdot k! \quad where \quad k=2 \nonumber \\
& = & \frac{(2+s)(2+s-1)\cdot \cdot \cdot (3)}{\rho ^{2}}\cdot \rho ^{k}k!(k=2) \nonumber \\
& = & \triangle _{0}\cdot \rho ^{2}\cdot 2! \quad  where \quad
\triangle _{0}=\frac{(2+s)(2+s-1)\cdot \cdot \cdot(3)}{\rho ^{2}}.
\nonumber
\end{eqnarray*}
It follows that $\int_{0}^{\infty }\frac{r^{k^{^{\prime
}}}}{e^{r}}dr\leq \triangle _{0}\cdot \rho ^{k}\cdot k!$. for all $k\geq 2$. \\
\\
\rm{(ii)}\em  Assume $\beta>0,\quad Then\quad $m$=\left\lceil
\frac{\beta }{2}\right\rceil$
 and $k\geq 2m+2.$ This gives
\[
\frac{k+s+1}{k+1}\leq \frac{2m+3+s}{2m+3}=1+\frac{s}{2m+3}.
\]
Let $\rho =1+\frac{s}{2m+3}$. Then
\[\int_{0}^{\infty }\frac{r^{k^{^{\prime \prime }}+1}}{e^{r}}dr\leq\triangle _{0}\cdot \rho ^{k+1}\cdot (k+1)!
\]
if $\int_{0}^{\infty }\frac{r^{k^{^{\prime \prime }}}}{e^{r}}dr\leq \triangle _{0}\cdot \rho ^{k}\cdot k!$. The smallest $k^{^{\prime \prime }}$ is $k_{0}^{^{\prime \prime }}$=2$m$+2+$s$ when $k$=2$m$+2. Now,\\
$\int_{0}^{\infty }\frac{r^{k_{0}^{^{\prime \prime }}}}{e^{r}}dr $
\begin{eqnarray*}
&=&k_{0}^{^{\prime \prime }}!=(2m+2+s)(2m+1+s)\cdot \cdot \cdot
(2m+3)(2m+2)! \\
&=&\frac{(2m+2+s)(2m+1+s)\cdot \cdot \cdot (2m+3)}{\rho ^{2m+2}}\cdot \rho
^{2m+2}\cdot (2m+2)! \\
&=&\triangle _{0}\cdot \rho ^{2m+2}\cdot (2m+2)! \quad where \quad \triangle _{0}=\frac{(2m+2+s)(2m+1+s)\cdot \cdot \cdot (2m+3)}{\rho ^{2m+2}}.
\end{eqnarray*}
It follows that \quad $\int_{0}^{\infty }\frac{r^{k^{^{\prime }}}}{e^{r}}dr\leq
\triangle _{0}\rho ^{k}k! \quad for \quad all \quad k\geq 2m+2.$\\
\\
\em \bf Case2.\rm \quad
Assume $k^{^{\prime \prime }}<k$. Let $k^{^{\prime \prime }}=k-s$ where $s>0$. Then
\[
\int_{0}^{\infty }\frac{r^{k^{^{\prime }}}}{e^{r}}dr\leq
\int_{0}^{\infty }\frac{r^{k^{^{\prime \prime }}}}{e^{r}}dr=k^{^{\prime
\prime }}!=(k-s)!=\frac{1}{k(k-1)\cdot \cdot \cdot (k-s+1)}\cdot k!
\]
and
\begin{eqnarray*}
\int_{0}^{\infty }\frac{r^{k^{^{\prime
}}+1}}{e^{r}}dr &\leq & \int_{0}^{\infty}\frac{r^{k^{^{\prime \prime }+1}}}{e^{r}}dr \\
& = &(k^{^{\prime \prime
}}+1)!=(k-s+1)!=\frac{1}{(k+1)k\cdot \cdot \cdot (k-s+2)}\cdot (k+1)!.
\end{eqnarray*}
Note that
\[
\left\{ \frac{1}{(k+1)k\cdot \cdot \cdot (k-s+2)}\left/ \frac{1}{k(k-1)\cdot
\cdot \cdot (k-s+1)}\right\}\right. =\frac{k(k-1)\cdot \cdot \cdot (k-s+1)}{(k+1)k\cdot \cdot \cdot (k-s+2)}=\frac{k-s+1}{k+1}.
\]
\rm{(i)} \em \quad Assume $\beta <0,$ then $m=0$ and $k\geq 2$. Since $k^{^{\prime \prime }}=k-s\geq 1$ holds for all $k\geq 2,$ it must be that $s=1.$
Thus
\[
\frac{k-s+1}{k+1}=1-\frac{s}{k+1}=1-\frac{1}{k+1}\leq 1 \quad for\quad all\quad k\geq 2.
\]
Let $ \rho =1.$ Then
\[
\int_{0}^{\infty }\frac{r^{k^{^{\prime \prime }}+1}}{e^{r}}dr\leq
\triangle _{0}\rho ^{k+1}(k+1)! \quad if \quad \int_{0}^{\infty }\frac{r^{k^{^{\prime
\prime }}}}{e^{r}}dr\leq \triangle _{0}\rho ^{k}k!.
\]
The smallest $k^{^{\prime \prime }\text{ }}$is $k_{0}^{^{\prime \prime }%
\text{ }}=k_{0}-s=2-s.$ Now,
\begin{eqnarray*}
\int_{0}^{\infty }\frac{r^{k_{0}^{^{\prime \prime }}}}{e^{r}}%
dr=k_{0}^{^{\prime \prime }\text{ }}!=(2-s)!=1! & =& 1 \\
& = &\frac{1}{2}k!\quad where\quad k=2 \\
& = &\frac{1}{2}\rho ^{k}k! \\
& = &\triangle _{0}\rho ^{k}k! \quad where \quad \triangle _{0}=\frac{1}{2}.
\end{eqnarray*}
It follows that $\int_{0}^{\infty }\frac{r^{k^{^{\prime }}}}{e^{r}}dr\leq
\triangle _{0}\rho ^{k}k!$ for all $k\geq 2.$\\
\rm{(ii)}\em Assume $\beta >0$ Then $m$=$\left\lceil \frac{\beta }{2}\right\rceil$
 and $k\geq 2m+2.$ This gives
\[
\frac{k-s+1}{k+1}=1-\frac{s}{k+1}\leq 1.
\]
Let $\rho =1$
Then $\int_{0}^{\infty }\frac{r^{k^{^{\prime \prime }}+1}}{e^{r}}dr\leq\triangle _{0}\cdot \rho ^{k+1}\cdot (k+1)!$ if $\int_{0}^{\infty }\frac{r^{k^{^{\prime \prime }}}}{e^{r}}dr\leq \triangle _{0}\cdot \rho ^{k}\cdot k!$.\\
The smallest $k$ is $k_{0}=2m+2$. Hence the smallest $k^{^{\prime \prime }}$ is $k_{0}^{^{\prime \prime }}$=$k_{0}-s=2m+2-s$. Now,
\begin{eqnarray*}
\int_{0}^{\infty }\frac{r^{k_{0}^{^{\prime \prime }}}}{e^{r}}dr
& = &k_{0}^{^{\prime \prime }}!=(2m+2-s)=(k_{0}-s)! \\
& = &\frac{1}{k_{0}(k_{0}-1)\cdot \cdot \cdot (k_{0}-s+1)}\cdot (k_{0})! \\
& = &\triangle _{0}\cdot \rho ^{k_{0}}k_{0}! \quad where \quad \triangle _{0}=\frac{1}{(2m+2)(2m+1)\cdot \cdot \cdot (2m-s+3)}.
\end{eqnarray*}
It follows that $\int_{0}^{\infty }\frac{r^{k^{^{\prime }}}}{e^{r}}dr\leq
\triangle _{0}\rho ^{k}k!$ for all $k\geq 2m+2.$\\
\\
\em \bf Case3.\rm \quad Assume $k^{^{\prime \prime }}=k$.\quad Then
\[
\int_{0}^{\infty }\frac{r^{k^{^{\prime }}}}{e^{r}}dr\leq
\int_{0}^{\infty }\frac{r^{k^{^{\prime \prime }}}}{e^{r}}dr=k!
\quad and \quad \int_{0}^{\infty }\frac{r^{k^{^{\prime
}}+1}}{e^{r}}dr\leq (k+1)!.
\]
Let $\rho =1.$ \quad Then
\begin{eqnarray*}
\int_{0}^{\infty }\frac{r^{k^{^{\prime }}}}{e^{r}}dr\leq \triangle
_{0}\rho ^{k}k! \quad for\quad all\quad k \quad where \quad
\triangle _{0}=1.
\end{eqnarray*}
The lemma is now an immediate result of the three cases. \hspace{6cm} $\Box$\\
\\
\em \bf Remark. \rm
For the convenience of the reader, we should express the
constants $\triangle _{0}$ and $\rho $ in a clear
form. It's easily shown that
\begin{enumerate}
\item[(a)] $k^{^{\prime \prime }}>k$ \quad if \quad and \quad only \quad if \quad $n-\beta >3$;
\item[(b)] $k^{^{\prime \prime }}<k$\quad if \quad and \quad only \quad if \quad $n-\beta \leq 1$;
\item[(c)] $k^{^{\prime \prime }}=k$\quad if \quad and \quad only \quad if \quad $1<n-\beta \leq 3$
\\
where \quad $k^{^{\prime \prime }}$ \quad and \quad $k$ \quad are \quad as \quad in \quad the \quad proof \quad of \quad the \quad lemma.
\end{enumerate}
We thus have the following situations.
\begin{enumerate}
\item[(a)] $n$-$\beta >3$.\quad Let $s=\left\lceil \frac{n-\beta
-3}{2}\right\rceil$
. \quad Then \\
\rm{(i)}\em if \quad $\beta <0$, $\rho =\frac{3+s}{3}$ and $\triangle _{0}=%
\frac{(2+s)(1+s)\cdot \cdot \cdot 3}{\rho ^{2}};$\\
\rm{(ii)}\em if \quad $\beta >0, \rho
=1+\frac{s}{2\left\lceil \frac{\beta }{2}\right\rceil +3}$ \quad
and \quad $\triangle _{0}=\frac{(2m+2+s)(2m+1+s)\cdot \cdot
\cdot (2m+3)}{\rho ^{2m+2}}$ \quad where \quad $m=\left\lceil \frac{\beta }{2}%
\right\rceil .$\\\em
\item[(b)] $n$-$\beta \leq 1$. Let $s=-\left\lceil \frac{n-\beta -3}{2}%
\right\rceil$ . Then \\
\rm{(i)}\em if \quad $\beta <0, \rho =1$ \quad and \quad $\triangle _{0}=\frac{1}{2};$\\
\rm{(ii)}\em if \quad $\beta >0, \rho =1$ \quad and \quad $\triangle _{0}=\frac{1}{%
(2m+2)(2m+1)\cdot \cdot \cdot (2m-s+3)}$; where $m=\left\lceil \frac{\beta }{%
2}\right\rceil .$\\\em
\item[(c)] $1<n-\beta \leq 3$. Then $\rho
=1 \quad and \quad \triangle _{0}=1.$
\end{enumerate}

  Before introducing our main theorem, we need the following lemma which is taken directly from $\left[ 6\right]$.
\begin{lemma}
Let $Q,Y$, and $\gamma _{n}$ be as in Lemma1.1. Then, given a
point $x$ in $Q$ ,there is a measure $\sigma$  supported on $Y$ such
that
\[
\int_{R^{n}}p(y)d\sigma (y)=p(x)
\]
for all $p$ in $\mathcal{P}_{k},$ \quad and
\[
\int_{R^{n}}d\left\vert \sigma \right\vert (y)\leq e^{2n\gamma _{n}(k+1)}.
\]
\end{lemma}

 Now we need another lemma.
\begin{lemma}
For any positive integer $k$,
\[
\frac{\sqrt{(2k)!}}{k!}\leq 2^{k}.
\]
\end{lemma}
\em \bf Proof.\rm \quad
This inequality holds for $k=1$ obviously. We proceed by
induction.
\begin{eqnarray*}
\frac{\sqrt{\left[ 2(k+1)\right]
!}}{(k+1)!}=\frac{\sqrt{(2k+2)!}}{k!(k+1)}=
\frac{\sqrt{(2k)!}}{k!}\cdot\frac{\sqrt{(2k+2)\cdot(2k+1)}}{k+1} \nonumber \\
\leq \frac{\sqrt{(2k)!}}{k!}\cdot
\frac{\sqrt{(2k+2)^{2}}}{k+1}\leq 2^{k}\cdot
\frac{(2k+2)}{k+1}=2^{k+1}.
\end{eqnarray*}
\hspace{14.6cm}\ \ \  $\Box$

 We can now enter the core of our theory. In the following theorem and its proof, we try to maintain the symbols of [6] whenever possible so that the reader can easily compare our results to those of Madych and Nelson. In fact the original form of this theorem is essentially Madych and Nelson's work. They offer the existence of the constant $\lambda$ and the exponential-type expression of the error bound. We include it here just to make this paper more readable. Our main contribution is the calculation of (i)the constant $\lambda$ which greatly depends on $\rho$, (ii)the constant $\delta_{0}$ and (iii)the coefficient to the left of $\lambda^{\frac{1}{\delta}}$ in (2.2) which greatly depends on the dimension n and the shape parameters $c $ and $\beta$.\\
\\
{\bf Theorem 2.4.}\  \em Suppose $h$ is defined as in $\left(1.4\right)$ and $m$ is its
order of conditional positive
definiteness. Let $\mu$ be its corresponding measure as in $\left( 1.5\right)$ .Then, given a positive number $b_{0}$, there are positive constants $\delta _{0}$ and $\lambda ,0<\lambda <1$, which depend on $b_{0}$ for which the following is true:\\
If $f\in \mathcal{C}_{h,m}$ and s is the h spline that
interpolates $f$ on a subset $X$ of $R^{n}$, then
\begin{align}
\left\vert f(x)-s(x)\right\vert \leq 2^{\frac{%
n+\beta +1}{4}}\cdot \pi ^{\frac{n+1}{4}}\cdot \sqrt{n\alpha _{n}}\cdot c^{%
\frac{\beta }{2}}\cdot \sqrt{\triangle _{0}}\cdot \lambda ^{\frac{1}{\delta }%
}\cdot \left\Vert f\right\Vert _{h} \tag{2.2}
\end{align}
holds for all $x$ in a cube $E$ provided that \rm{(i)}\em $E$ has
side $b$ and $b\geq b_{0}$,\rm{(ii)}\em $0<\delta \leq \delta
_{0}$ and \rm{(iii)}\em every subcube of $E$ of side $\delta $
contains a point of $X.$
Here, $\alpha _{n}$ denotes the volume of the unit ball in $R^{n}$ and $c,$ $%
\triangle _{0}$ are as in $\left( 2.1\right) .$\\
\indent The numbers $\delta _{0}$ and $\lambda $ can be expressed specifically as
\[
\delta _{0}=\frac{1}{3C\gamma _{n}(m+1)},~\lambda =\left( \frac{2}{3}\right) ^{\frac{1}{%
3C\gamma _{n}}}
\]
where
\[
C=\max \left\{ 2\rho ^{\prime }\sqrt{n}e^{2n\gamma _{n}},~\frac{2}{3b_{0}}%
\right\} ,~\rho ^{^{\prime }}=\frac{\rho }{c}.
\]
The number $\rho $ can be found in the remark immediately following
Lemma2.1, and $\gamma _{n}$ was defined in Lemma1.1.\\
\\
\em \bf Proof.\rm \quad
Let $\rho$ and $\gamma _{n}$ be as in the statement of the theorem. Fix the parameter $c$ in (1.4). For any $b_{0}>0$, let
\[
B=2\rho ^{^{\prime }}\sqrt{n}e^{2n\gamma _{n}} \quad and \quad C=\max \left\{ B,~%
\frac{2}{3b_{0}}\right\} \quad where \quad \rho^{^{\prime }}=\frac{\rho }{c}.
\]
 Let $\delta_{0}$ be defined as in the statement of the theorem.

 We start our proof with a crucial inequality which is a result of Theorem4.2 of [5]. Let $E$ be the cube mentioned in the theorem. For any $x\in E$,
\begin{align}
\left\vert f(x)-s(x)\right\vert \leq c_{k}\left\Vert
f\right\Vert _{h}\int_{R^{n}}\left\vert y-x\right\vert ^{k}d\left\vert
\sigma \right\vert (y) \tag{2.3}
\end{align}
whenever $k>m$, where $\sigma $ is any measure supported on $X$ such that\\
\begin{align}
\int_{R^{n}}p(y)d\sigma(y)=p(x)\tag{2.4}
\end{align}
for all polynomials $p$ in $\mathcal{P}_{k-1}.$ Here
\[
c_{k}=\left\{ \int_{R^{n}}\frac{\left\vert \xi \right\vert ^{2k}}{(k!)^{2}}%
d\mu (\xi )\right\} ^{\frac{1}{2}}
\]
whenever \quad $k>m.$ \quad By $\left( 2.1\right) ,$ \quad for all
\quad $2k\geq 2m+2,$
\begin{enumerate}
\item[(2.5)]
\quad \quad \quad $c_{k}=\left\{ \int_{R^{n}}\frac{\left\vert \xi \right\vert
^{2k}}{(k!)^{2}}d\mu (\xi )\right\} ^{\frac{1}{2}}$
\end{enumerate}
\begin{align}
\quad \quad \quad \leq \frac{1}{k!}\cdot 2^{\frac{n+\beta +1}{4}}\cdot \pi ^{\frac{n+1}{4}}\cdot \sqrt{n\alpha _{n}}\cdot c^{\frac{\beta -2k}{2}}\cdot \sqrt{\triangle _{0}}\cdot \rho ^{k}\cdot
\sqrt{(2k)!} \nonumber
\end{align}
\begin{align}
\leq 2^{\frac{n+\beta +1}{4}}\cdot \pi ^{\frac{n+1}{4}}\cdot \sqrt{n\alpha
_{n}}\cdot c^{\frac{\beta }{2}}\cdot c^{-k}\cdot \sqrt{\triangle _{0}}\cdot
(2\rho )^{k} \nonumber
\end{align}
due to Lemma2.3.

 In order to develop the inequality (2.3) into (2.2), we have to find a bound for the value
\[
I:=c_{k}\int_{R^{n}}\left\vert y-x\right\vert^{k} d\left\vert \sigma \right\vert
(y).
\]
For this we have to appeal to Madych and Nelson's theory in [5] an [6].

 Let $\delta>0 $ be as in the statement of the theorem. Since $%
\delta \leq \delta _{0}$, one easily finds that  $0<3Cr _{n}\delta \leq \frac{1}{m+1}$. Obviously we can choose an integer $k\geq m+1$ so that
\[
1\leq 3C\gamma _{n}k\delta \leq 2.
\]
By the definition od $C$ and simple calculation we get $\gamma _{n}k\delta \leq b_{0}$ for such a $k.$ Let $Q$ be any
cube which contains $x,$ has side $\gamma _{n}k\delta ,$ and is contained in
$E.$ Subdivide $Q$ into $\left( \gamma _{n}k\right) ^{n}$ subcubes
of side $\delta$. By hypothesis each of these subcubes must contain a point of $X$.
Select arbitrarily a point of $X$ from each such subcube and let $Y$ denote the set of these points. As a result of Lemma2.2, there exists a measure $\sigma $ supported on $Y$ satisfying $\left( 2.4\right) $ and the following inequality.
\begin{align}
\int_{R^{n}}d\left\vert \sigma \right\vert (y)\leq e^{2n\gamma
_{n}k}.\tag{2.6}
\end{align}
With the help of this measure and a bound on $I$ can be constructed.

 Because the support of
$\sigma $ is contained in $Q$ whose diameter is $\sqrt{n}\gamma _{n}k\delta
, $ by (2.5) and (2.6), we get
\begin{enumerate}
\item[(2.7)]
$ \quad \quad \quad I\leq 2^{\frac{n+\beta +1}{4}}\cdot \pi ^{\frac{n+1}{4}}\cdot
\sqrt{n\alpha _{n}}\cdot c^{\frac{\beta }{2}}\cdot c^{-k}\cdot \sqrt{%
\triangle _{0}}\cdot (2\rho )^{k}(\sqrt{n}\gamma _{n}k\delta
)^{k}e^{2n\gamma _{n}k}$
\end{enumerate}
\begin{align}
 \leq \left( C\gamma _{n}k\delta \right) ^{k}(2^{\frac{n+\beta +1}{4}}\cdot
\pi ^{\frac{n+1}{4}}\cdot \sqrt{n\alpha _{n}}\cdot c^{\frac{\beta }{2}}\cdot
\sqrt{\triangle _{0}}).\nonumber
\end{align}
Since
\[
C\gamma _{n}k\delta \leq \frac{2}{3} \quad and \quad k\geq \frac{1}{3C\gamma
_{n}\delta },
\]
(2.7) implies that
\[
I\leq \left[ (\frac{2}{3})^{\frac{1}{3C\gamma _{n}}}\right] ^{\frac{1}{%
\delta }}\cdot (2^{\frac{n+\beta +1}{4}}\cdot \pi ^{\frac{n+1}{4}}\cdot
\sqrt{n\alpha _{n}}\cdot c^{\frac{\beta }{2}}\cdot \sqrt{\triangle _{0}}).
\]
This together with (2.3) gives
\[
\left\vert f(x)-s(x)\right\vert \leq 2^{\frac{n+\beta +1}{4}}\cdot \pi ^{%
\frac{n+1}{4}}\cdot \sqrt{n\alpha _{n}}\cdot c^{\frac{\beta }{2}}\cdot \sqrt{%
\triangle _{0}}\cdot \lambda ^{\frac{1}{\delta }}\cdot \left\Vert
f\right\Vert _{h},
\]
where
\[
\lambda =\left( \frac{2}{3}\right) ^{\frac{1}{3C\gamma _{n}}}.
\] $\hspace{14.75cm}\ \ \Box$ \\
\\
{\bf Remark}. The value of $\lambda$ can now easily be obtained by its very definition in Theorem2.4. For example, if we fix $n=1,b_{0}=1,\beta=1$ and let $c=1,5,$ and $10$, then $\lambda$ will be 0.999381, 0.99691 and 0.99383, repectively. This shows an important fact that the crucial constant $\lambda$ has a close relationship with $c$ and a criterion for the optimal choice of c may be developed.\\
\indent What's noteworthy is that in Theorem2.4 the parameter $\delta$ is not the
generally used fill distance. For easy use we should transform the theorem
into a statement described by the fill distance.

 Let
\[
d\left( \Omega ,X\right) =\sup_{y\in \Omega }\,\inf_{x\in X}\left\vert y-x\right\vert
\]
be the fill distance. Observe that every cube of side $\delta $ contains a ball of radius $\frac{\delta }{2}.$
Thus the subcube condition in Theorem2.4 is satisfied when $\delta
=2d(E,X).$ More generally, we can easily conclude the following:\\
\\
{\bf Corollary2.5}\ \   \em Suppose $h$ is defined as in $\left( 1.4\right) $ and $m$ is its order of conditional positive definiteness. Let $\mu $ be its corresponding measure
as in $\left( 1.5\right) .$ Then, given a positive number $b_{0},$%
 there are positive constants $d_{0}$ and $\lambda^{^{\prime }},$ $0<\lambda ^{^{\prime }}<1,$ which depend on $b_{0}$ for which the following is true: If $f\in \mathcal C_{h,m}$ 
 and $s$ is the $h$ spline that interpolates $f$ on a subset $X$ of $R^{n},$ then
\begin{align}
\left\vert f(x)-s(x)\right\vert \leq 2^{\frac{n+\beta +1}{4}%
}\cdot \pi ^{\frac{n+1}{4}}\cdot \sqrt{n\alpha _{n}}\cdot c^{\frac{\beta }{2}%
}\cdot \sqrt{\triangle _{0}}\cdot (\lambda ^{^{\prime }})^{\frac{1}{d}}\cdot
\left\Vert f\right\Vert _{h}\tag{2.8}
\end{align}
holds for all $x$ in a cube $E\subseteq \Omega ,$ where $\Omega $ is a set
which can be expressed as the union of rotations and translations of a fixed
cube of side $b_{0},$ provided that \rm{(i)}\em $E$ has side $b\geq
b_{0},~$\rm{(ii)}\em $0<d\leq$ $d_{0}$ and \rm{(iii)}\em  every subcube of $E$ of side $2d$
contains a point of $X.$ Here, $\alpha _{n}$ denotes the volume of the unit
ball in $R^{n}$ and $c,~\triangle _{0}$ are as in $\left( 2.1\right) .$
Moreover $d_{0}=\frac{\delta _{0}}{2}$ and $\lambda ^{^{\prime }}=\sqrt{%
\lambda }$ where $\delta _{0}$ and $\lambda $ are as in Theorem2.4.\\
\\
 {\bf Proof.}  \rm  Let $d_{0}=\frac{\delta _{0}}{2}$ and $\delta =2d.$ Then $0<d\leq d_{0}$ iff
$0<\delta \leq \delta _{0}.$ Our corollary follows immediately by noting that
$\lambda ^{\frac{1}{\delta }}=\lambda ^{\frac{1}{2d}}=\sqrt{\lambda }^{\frac{%
1}{d}}=(\lambda ^{\prime })^{\frac{1}{d}}.$\ \ \ \   $\hspace{8.1cm} \Box$ \\
\\
\em \bf Remark. \rm
(a)The space $\mathcal C_{h,m}$ probably is unfamiliar
to most people. It was introduced by Madych and Nelson in $\left[ 4\right] $%
 and $\left[ 5\right] .$ Later Luh made characterizations for it in $\left[ 1\right] $ and $\left[ 2\right] .$
Many people think that it is defined by Gelfand and Shilov's definition of
generalized Fourier transform, and is therefore difficult to deal with. This
is not true. In fact, it can be characterized by Schwartz's definition of
generalized Fourier transform. The situation is not so bad. Moreover, many
people mistake $\mathcal C_{h,m}$ to be the closure of Wu and Schaback's
function space defined in $\left[ 9\right] .$ This is also
not true. The two spaces have very subtle connection. Luh has also made a
clarification for this problem. For further details, please see $\left[ 3%
\right] $. 
(b)In the proof of Lemma2.1 although the integration was obtained by approximation in the fifth line, the gap for $\xi$ near the origin can be made arbitrarily small by decreasing $\delta_{0}$ of Theorem2.4. Note that the integrand collapses to zero for $0<|\xi|<1$ as $k\rightarrow \infty$. The third paragraph of page 9 shows that $k\rightarrow \infty$ as $\delta\rightarrow 0$. Hence it's harmless. As for $|\xi|\geq 1$, this approximation is a commonly used approach and the proof of Theorem2.4 shows that the gap can be ignored from the viewpoint of error bound.
(c)Seemingly our main results Theorem2.4 and Corollary2.5 are quite complicated. However they strongly promise the birth of a useful set of criteria for the optimal choice of the shape parameters c and $\beta$ contained in the multiquadrics. All these will be seen in the forthcoming papers of the author.
\begin{center}
REFERENCES
\end{center}

\begin{enumerate}
\item Lin-Tian Luh, The Equivalence Theory of Native Spaces,
Approx. Theory and its Applications, 2001, 17:1, 76-96.

\item Lin-Tian Luh, The Embedding Theory of Native Spaces, Approx.
Theory and its Applications, 2001, 17:4, 90-104.

\item Lin-Tian Luh, On Wu and Schaback's Error Bound, Inter. J. Numer. Methods and Applications, Vol.1, N0.2, 2009, 155-174.

\item W.R.Madych and S.A.Nelson, Multivariate interpolation and
conditionally positive definite function, Approx. Theory Appl. 4, No.
4, 1988, 77-89.

\item W.R.Madych and S.A.Nelson, Multivariate interpolation and
conditionally positive definite function, II, Math. Comp. 54, 1990, 211-230.

\item W.R.Madych and S.A.Nelson, Bounds on Multivariate Polynomials and
Exponential Error Estimates for Multiquadric Interpolation, J. Approx.
Theory 70, 1992, 94-114.

\item R. Schaback and H. Wendland, Characterization and Construction of
Radial Basis Functions, preprint.

\item H. Wendland,
Scattered Data Approximation, Cambridge University Press, 2005.

\item Z. Wu and R. Schaback, Local Error Estimates for Radial Basis Function
Interpolation of Scattered Data, IMA J. of Numerical Analysis,13, 1993,
13-27.
\end{enumerate}

\parbox[t]{3.5in}{Lin-Tian Luh\\ Department of Mathematics\\ Providence University\\ 
  Shalu Area\\ Taichung City\\ Taiwan\\ltluh@pu.edu.tw\\ }

\end{document}